\documentclass [12pt]{amsart}

\usepackage{amssymb,amsxtra,amsfonts}
\usepackage{epsf}              %
\usepackage{epsfig}             %
\usepackage{graphics}
\def\reE@DeclareMathSymbol#1#2#3#4{%
    \let#1=\undefined
    \DeclareMathSymbol{#1}{#2}{#3}{#4}}
\DeclareSymbolFont{symbolsC}{U}{txsyc}{m}{n}
\SetSymbolFont{symbolsC}{bold}{U}{txsyc}{bx}{n}
\DeclareFontSubstitution{U}{txsyc}{m}{n}
\reE@DeclareMathSymbol{\strictiff}{\mathrel}{symbolsC}{76}

\usepackage{verbatim}  %for \comment
\usepackage{bbding}   % for \bigoasterisk----now use made-up version 

\def\reE@DeclareMathSymbol#1#2#3#4{%
    \let#1=\undefined
    \DeclareMathSymbol{#1}{#2}{#3}{#4}}
\DeclareSymbolFont{symbolsC}{U}{txsyc}{m}{n}
\SetSymbolFont{symbolsC}{bold}{U}{txsyc}{bx}{n}
\DeclareFontSubstitution{U}{txsyc}{m}{n}
\reE@DeclareMathSymbol{\strictiff}{\mathrel}{symbolsC}{76}

\newcommand\beq{\begin{equation}}
\newcommand\eeq{\end{equation}}
\newcommand\bal{\begin{align*}}
\newcommand\eal{\end{align*}}   %why does this not work??
\newcommand\bmx{\left(\begin{matrix}}
\newcommand\emx{\end{matrix}\right)}
\newcommand\bsmx{\left(\begin{smallmatrix}}
\newcommand\esmx{\end{smallmatrix}\right)}

\newcommand{\spq}{/\!\!/}

\newcommand{\st}{\ \bigl\vert\ }

\providecommand{\im}{\text{\rm Im}}

\def\part#1{\frac{\partial\phantom{q}}{\partial#1}}

\newcommand{\glue}[1]{\underset{#1}{\strictiff}}
\newcommand{\fus}{\circledast}

\newcommand{\HH}{\text{\rm H}}

\newcommand{\Id}{\text{\rm Id}}

\newcommand{\Lie}{{\mathop{\rm Lie}}}
\newcommand{\Vect}{{\mathop{\rm Vect}}}             % see \vrt for Vert 
\newcommand{\pap}[2]{{\ _{#1}\cA_{#2}}}

\newcommand{\gah}{\pap{G}{H}}  %{\ _G\cA_H}}
\newcommand{\Ad}{{\mathop{\rm Ad}}}
\newcommand{\ad}{{\mathop{\rm ad}}}

\DeclareMathOperator{\Hom}{Hom}         % this looks to be the correct way to do this (put a * ``...otar*{'' if want things underneath

\newcommand{\GL}{{\mathop{\rm GL}}}

\newcommand{\U}{{\rm {U}}}	%bizarre-- mathop shifts it down!
\newcommand{\irr}{{\rm irr}}

\renewcommand{\ker}{\mathop{\rm Ker}}
\newcommand{\bs}{{\bf S}}

\newcommand{\IC}{\mathbb{C}}

\newcommand{\A}{\mathcal{A}}
\newcommand{\cA}{\mathcal{A}}
\newcommand{\cB}{\mathcal{B}}
\newcommand{\cC}{\mathcal{C}}

\newcommand{\ch}{\eta}     %{\hslash}    %\hbar, hslash..mathcal fails

\newcommand{\cM}{\mathcal{M}}

\newcommand{\cP}{\mathcal{P}}
\newcommand{\cp}{\mathcal{P}}

\newcommand{\cR}{\mathcal{R}}
\newcommand{\cS}{\mathcal{S}}

\newcommand{\cU}{\mathcal{U}}
\newcommand{\cu}{\mathcal{U}}

\renewcommand{\u}{       \mathfrak{u}     }
\newcommand{\g}{       \mathfrak{g}     }

\newcommand{\lt}{\mathfrak{t}}

\newcommand{\so}{       \mathfrak{so}     }

\renewcommand{\sp}{{\mathfrak{sp}}}

\newcommand{\h}{\mathfrak{h}}

\newcommand{\lu}{\mathfrak{u}}

\newcommand{\wt}{\widetilde}

\newcommand{\al}{\alpha}

\newcommand{\be}{\beta}
\newcommand{\ga}{\gamma}
\newcommand{\de}{\delta}

\newcommand {\eps}{\varepsilon}

\newcommand{\La}{\Lambda}

\newcommand{\ze}{\zeta}

\renewcommand{\bar}{\overline}

\makeatletter
 \newlength{\typesize}
\makeatother

\newlength{\vvoff}
\newlength{\hhoff}

\def\underset#1#2{\ \smash{\mathop{ #2 }\limits_{#1}}\ }

\newcommand{\pf}{\begin{bpf}}

\newcommand{\pfms}{\begin{bpfms}}
\newcommand{\epf}{\end{bpf}\hfill$\square$\\}           % end proof
\newcommand{\epfms}{\end{bpfms}\hfill$\square$\\}       % end proof

\newcommand{\idea}{\begin{bidea}}

\newcommand{\eidea}{\end{bidea}\hfill$\square$\\}           % end proof

\newcommand{\sk}{\begin{bsk}}    %type: \sk ..... \esk

\newcommand{\esk}{\end{bsk}\hfill$\square$\\}           % end sketch
\newcommand{\sketch}{\begin{bsketch}}%type: \sketch ..... \esketch

\newcommand{\esketch}{\end{bsketch}\hfill$\square$\\}

\theoremstyle{plain}  \newtheorem{hypo}{Hypothesis}[section]
\newtheorem{thm}[hypo]{Theorem}

\newtheorem {defn}[hypo]{Definition}

\theoremstyle{definition} \newtheorem{rmk}[hypo]{Remark}

\newtheorem{eg}[hypo]{Example}

\begin{document}
\title[Through the analytic halo]{Through the analytic halo:\\ Fission via irregular singularities}
\author{Philip Boalch}

\dedicatory{To Bernard Malgrange on his 80th birthday}

\maketitle

\tableofcontents

\nocite{DMR-ci}

\section{Introduction}

The aim of the article is to tentatively initiate a direction of research concerning connections on bundles on Riemann surfaces, where the structure group of the bundle may vary in different regions of the surface. 
In the simplest case there will be certain loops (``halos") drawn on the surface across which the structure group will be broken to a subgroup\footnote{
In fact this situation arises quite naturally by thinking about irregular connections in a slightly novel way (see Remark  \ref{rmk: irreg}), although for the most part we will work more topologically.}. %

Here we will describe such moduli spaces as complex symplectic manifolds, extending the topological/Betti viewpoint---the spaces to be described here generalise the character varieties of Riemann surfaces. We will postpone until later further investigation of finer properties, such as the existence of hyperk\"ahler metrics, correspondence with Higgs bundles, and the possibility of extending the geometric Langlands program (we are, after all, generalising the spaces on the Galois side of this correspondence).

Philosophically we are interested in trying to generalise the first nonabelian cohomology 
$$\cM = \HH^1(\Sigma, G)$$
of a smooth projective curve $\Sigma$ with coefficients in a complex reductive group $G$. 
By work of Hitchin, Donaldson, Corlette, Simpson and others
such a cohomology space may be realised in various ways with different algebraic structures. The rough picture is as follows. The {\em Betti realisation}
is as the space
\beq \label{eq: hom p1 g}
\Hom(\pi_1(\Sigma),G)/G
\eeq
of conjugacy classes of representations of the fundamental group of $\Sigma$.
By the Riemann--Hilbert correspondence
this is isomorphic to a space of holomorphic connections on $G$-bundles on $\Sigma$
(the {\em de\! Rham realisation}).
By the nonabelian Hodge theorem, this is also a space of 
Higgs fields on $G$-bundles on $\Sigma$ (the {\em Dolbeault realisation}).
The different complex structures on the space $\cM$ may be expressed in terms of the existence of a natural hyperk\"ahler metric on $\cM$
(see \cite{Hit-sde, simp-hfnc}).

Simpson \cite{Sim-hboncc} extended the nonabelian Hodge correspondence to the case of a punctured curve; essentially one is now considering representations of the fundamental group of the punctured curve (or connections/Higgs fields with simple poles). This may still be understood in terms of hyperk\"ahler metrics, considered for example by Nakajima \cite{Nak}, but one needs to first restrict to fundamental group representations taking loops around each puncture into fixed  conjugacy classes of $G$.
Said differently
$$\Hom(\pi_1(%
{\Sigma\setminus\{\text{$m$ points}\}}),G)/G$$
has a natural holomorphic Poisson structure and its (generic) symplectic leaves,
which are obtained by fixing the local monodromy conjugacy classes, are hyperk\"ahler.

It turns out that one may replace the fundamental group here by the 
{\em wild fundamental group}
of Martinet--Ramis \cite{MR91}
(abstractly this is just the Tannaka group of the Tannakian category of meromorphic connections on vector bundles on $\Sigma$).
Namely, roughly speaking,
$$\Hom(\pi^s_1(\Sigma),G)/G$$
has a natural holomorphic Poisson structure where $\pi^s_1$ is the wild fundamental group. 
The symplectic leaves of this Poisson structure are finite dimensional and correspond to meromorphic connections with fixed formal type at each pole 
(naturally extending the notion of fixing the local monodromy conjugacy class above).
In \cite{wnabh} it was shown that the (sufficiently generic) symplectic leaves have 
hyperk\"ahler metrics and a correspondence between meromorphic connections and Higgs fields was established; %
 a ``wild'' nonabelian Hodge correspondence on curves.

Thus in brief there are various modifications of the first slot in 
\eqref{eq: hom p1 g} that may be made.
Now we would like to try to modify the second slot.
As mentioned above the picture in mind is of a Riemann surface with certain loops (``halos") drawn on it dividing it into pieces.
On each piece one is free to choose a complex reductive group, provided that---in the simplest case---the group on one side of each halo should be the stabiliser of a semisimple element of the Lie algebra of the group on the other side.
Our goal here is to construct such spaces as complex symplectic manifolds.

\begin{rmk}
The term ``analytic halo'' is borrowed from Martinet--Ramis \cite{MR91}. See also 
Deligne--Malgrange--Ramis \cite{DMR-ci}.
\end{rmk}

\section{Strategy}

Symplectic spaces of flat connections on $G$-bundles on surfaces have been intensively studied in recent decades and there are many different approaches. The finite dimensional ``quasi-Hamiltonian''  approach of Alekseev--Malkin--Meinrenken \cite{AMM} involves fusing together some basic pieces and then performing a reduction to obtain the symplectic moduli space. This motivated their theory of Lie group valued moment maps.
Some familiarity with this theory will aid the reader (see \cite{AMM}, and  \cite{saqh} for the holomorphic version).
We view it as a convenient algebraic framework to make precise analogues of various analytic operations involving loop groups.

Given a compact Riemann surface $\Sigma$ with one boundary circle one obtains a quasi-Hamiltonian $G$-space by taking the moduli space of flat connections on $G$-bundles on $\Sigma$ with a framing at one point of the boundary. 
Similarly a surface with $m$ boundary components leads to a quasi-Hamiltonian $G^m$-space, by including a framing at one point in each boundary component.

The fusion operation corresponds to gluing two surfaces with one boundary component onto two of the boundary circles of a three-holed sphere (and thus obtaining a new quasi-Hamiltonian $G$-space corresponding to the resulting surface---which still has just one boundary component).
This puts a ring structure on the category of quasi-Hamiltonian $G$-spaces, the identity for which is the space corresponding to the disk.
Fusing with a disk corresponds to gluing on an annulus to the original surface, which does not change the moduli space of flat connections.
Now the annulus has two boundary components so corresponds naturally to a quasi-Hamiltonian $G\times G$-space (the `double' of \cite{AMM}), whose moment map may be written as follows:
$$ D(G) = G\times G; \quad \mu(C,h) = (C^{-1}hC,h^{-1})\in G\times G.$$
The two components of the moment map correspond to the holonomies of connections on the annulus around the two boundary components.

Now imagine a surface with a halo on it (i.e. an embedded circle), and consider the annulus given by a tubular neighbourhood of the halo. If the structure group is $G$ on one side of the halo and $H$ on the other side then the holonomies around the boundary components will be in $G$ and $H$ respectively. Thus there should  be a quasi-Hamiltonian $G\times H$-space attached to the annulus 
with moment map given by the two holonomies as before, and containing something extra  (a little more complicated) related to crossing the halo, going from one boundary component to the other.

Suppose (hypothetically) that such a quasi-Hamiltonian $G\times H$-space $\gah$ exists.  
Then we obtain a symplectic manifold by taking an arbitrary quasi-Hamiltonian $G$-space $M_G$ and quasi-Hamiltonian $H$-space $M_H$ and gluing them to the annulus, i.e. by performing the fusions
$$M_G \underset{G}{\fus} \gah \underset{H}{\fus} M_H$$ 
(to obtain a new quasi-Hamiltonian $G\times H$-space) and 
then reducing by $G\times H$ to obtain a symplectic manifold (if it is a manifold). 

\begin{rmk}
Note that the operation
$$M_G\mapsto (M_G \underset{G}{\fus} \gah)\spq G$$
will associate a quasi-Hamiltonian $H$-space to any quasi-Hamiltonian $G$-space $M_G$,
and similarly one may obtain a quasi-Hamiltonian $G$-space from a quasi-Hamiltonian $H$-space.
\end{rmk}

In this way a general surface with some (nonintersecting) halos drawn on it and chosen structure groups leads to a symplectic manifold by cutting it up into pieces and gluing as above. 
Thus the general problem reduces to that of establishing the existence of the quasi-Hamiltonian $G\times H$-spaces $\gah$ attached to an annulus containing a halo.
This will be established in the following section.

In the section after next we will consider the case when $H$ is a product of groups;  one may then glue on a quasi-Hamiltonian $H_i$ space for each factor $H_i$ of the group $H$.
As will be explained, this leads to a {\em fission} picture, breaking $G$ in to the pieces $H_i$; 
The space $\gah$ should perhaps best be pictured in terms of a surface having  a boundary circle  for each factor of $H$ and also for $G$.

\section{The Fission Spaces}

Let $G$ be a connected complex reductive group with Lie algebra $\g$.
Choose a nondegenerate invariant bilinear form $(\ ,\ )$ on $\g$. 
Let $A\in \g$ be a semisimple element, and let $H\subset G$ be the stabiliser of $A$ under the adjoint action.
(Then $H$ is again a connected complex reductive group.)

Suppose we are given a Cartan subalgebra $\lt\subset \g$ containing $A$
and a choice of positive roots $\cR_+\subset \cR\subset \lt^*$, where $\cR$ is the set of roots of $\g$. Write $\cR_-=-\cR_+$ for the corresponding negative roots.

Since $A$ is semisimple one has a vector space direct sum:
$$\g = \im(\ad_A)\oplus\ker(\ad_A)$$
where $\ker(\ad_A)$ is the Lie algebra $\h$ of $H$.
The complementary subspace  $\im(\ad_A)$ is stabilised by $\lt$ and breaks up as a direct sum of the (one-dimensional) root spaces of $\g$ that it contains, and so we may write
$$\im(\ad_A) = \lu_+\oplus \lu_-$$
for the subspaces corresponding to positive and negative roots.
Now $\lu_+$ and $\lu_-$ are nilpotent Lie subalgebras of $\g$ and we may exponentiate them to 
obtain unipotent subgroups $U_+, U_-\subset G$. 
The desired quasi-Hamiltonian $G\times H$-space is as follows.

Let $\theta,\overline \theta\in\Omega^1(G,\g)$
denote the left and right invariant Maurer--Cartan forms on $G$ respectively, and let 
$\theta_\pm,\overline \theta_\pm\in\Omega^1(U_\pm,\u_\pm)$
denote the Maurer--Cartan forms on $\U_\pm$.

\begin{thm} \label{thm: qh}
The space
$$\gah := G\times U_-\times U_+\times H$$
is a complex quasi-Hamiltonian $G\times H$-space, 
with $G\times H$ action:
$$(g,k)\cdot (C,\ u_-,\ u_+,\ h)  = 
(kC g^{-1},\  ku_-k^{-1},\  ku_+k^{-1},\  khk^{-1}),$$
(where $(g,k)\in G\times H$), with moment map:
$$\mu(C,\ u_-,\ u_+,\ h) = 
(C^{-1} p C,\  h^{-1}) \in G\times H,$$
where  $p=u_-^{-1}hu_+\in G$, 
and with holomorphic two-form $\omega$ given by:

$$ \omega =\frac{1}{2} 
\left(\overline\gamma, \Ad_p\overline \gamma\right)
+\frac{1}{2} \left(\overline \cu_-, \Ad_h\overline \cu_+\right)
+\frac{1}{2} \left( \overline \gamma,\cp+\overline \cp\right)
$$
where
$\overline \gamma  = C^*(\overline\theta)$,
$\cp  = p^*(\theta),\bar\cp  = p^*(\bar\theta)$, and 
$\bar\cU_\pm = u_\pm^*(\bar\theta_\pm)$.
\end{thm}
\pf This may be verified directly (see the appendix). \epf

\begin{eg}
Suppose $A=0$. Then $H=G$ and both  of $U_+$ and $U_-$ are a point. 
Then $\gah = G\times G $ is the double, which is one of the basic examples of Alekseev--Malkin--Meinrenken \cite{AMM}.
\end{eg}
\begin{eg}
Suppose $A$ is regular. 
Then $H$ is a maximal torus of $G$ and $U_\pm$ are the unipotent radicals of a pair of opposite Borels in $G$.
Now (by definition) as a space 
the standard dual Poisson Lie group $G^*$ 
of $G$ is a covering of $U_+\times U_-\times H$ and one may pull back the quasi-Hamiltonian structure on $\gah$ to obtain the structure of  quasi-Hamiltonian 
$G\times H$-space on the product  $G\times G^*$. 
This is one of the basic examples of \cite{saqh}.
\end{eg}

Thus, for more general $H$, the spaces $\gah$ interpolate between the above two examples. They always have dimension equal to twice that of $G$.

\section{Fission}

Now suppose that  $H=H_1\times H_2$ is written as a product of two groups.
(The generalisation to arbitrarily many factors is immediate.)

Given a reductive group $G$ and a surface $\Sigma$ with one boundary circle and with a marked point $x$ on the boundary, let 
$$M_G(\Sigma)$$ 
denote the quasi-Hamiltonian $G$-space obtained as usual by taking the space of flat $G$-connections on $\Sigma$ with a framing at $x$.

Note that if $H=H_1\times H_2$ then 
$$M_H(\Sigma) \cong M_{H_1}(\Sigma)\times M_{H_2}(\Sigma)$$ 
since specifying a $H$-connection on $\Sigma$ is the same as specifying a pair consisting of a $H_1$-connection on $\Sigma$ and a $H_2$-connection on $\Sigma$.

Note also that in general the product $M_{H_1}\times M_{H_2}$ is a quasi-Hamiltonian $H$-space for any quasi-Hamiltonian $H_i$-spaces $M_{H_i}$. In particular we may take any two surfaces $\Sigma_1,\Sigma_2$ with one boundary component and set
$$M_H = M_{H_1}(\Sigma_1)\times M_{H_2}(\Sigma_2).$$ 

Thus if $M_G = M_G(\Sigma)$ then the  space

$$M_G \underset{G}{\fus} \gah \underset{H}{\fus} M_H\quad$$
is obtained by gluing together the three surfaces $\Sigma, \Sigma_1, \Sigma_2$.
This suggests that, rather than thinking of $\gah$ in terms of an annulus, one should think of a surface with three boundary components labelled by $G,H_1,H_2$ respectively.
We think of this surface as 
the product $Y\times S^1$ of a $Y$-shaped piece and a circle.
The surfaces $\Sigma_1$ and $\Sigma_2$ do not interact---one is in effect sewing {\em both} of them on to the boundary of $\Sigma$.

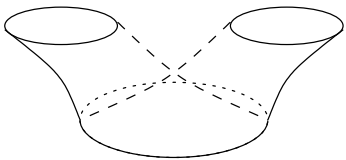
\begin{figure}[ht]
	\centering
	\input{fission2.pstex_t}
\end{figure}

\begin{rmk}
In the case when $G=H=H_1\times H_2$ (i.e. $G$ is itself a product and we take $A=0$) the global picture reduces to that of two surfaces, each with flat $H_i$-connections respectively, and with some regions of the two surfaces identified (where one thinks in terms of a single $G$-connection rather than two connections). The general case does not decouple in this way.
\end{rmk}

Thus if $H=\prod_1^n H_i$ the quasi-Hamiltonian $G\times H$-space $\gah$ breaks the symmetry group from $G$ into the pieces $H_i$. We view this as a `fission' operation,
complementary (and not inverse) to the usual fusion operation.
Rather, fission breaks the group from $G$ to $H=\Pi H_i$, and this then 
yields the opportunity to fuse with quasi-Hamiltonian $H_i$-spaces.
(This was referred to as  ``fusion on the other side of the analytic halo''  in 
\cite{rsode} footnote 3.)

For example note that 
taking $\Sigma_1$ to be a disc gives a way to kill a factor of $H$, reducing $G$ to a proper subgroup of $H$.

Note also 
that there are many pairs of reductive groups having isomorphic subgroups $H$, so may be glued together, possibly via intermediate surfaces. E.g. in the simplest case any two groups of the same rank, taking $H$ to be a maximal torus, although it is easy to construct examples with nonabelian $H$.\footnote{
The stabiliser of a semisimple Lie algebra element is a Levi factor of a parabolic subgroup, and these correspond to subsets of the nodes of the Dynkin diagram of $G$.
E.g. one could  use the fact that the Dynkin diagrams $B_2$ and $C_2$ coincide to glue groups having Lie algebras 
$\so_{2n+1}(\IC)$ and $\sp_{2n}(\IC)$.}

\section{Gluing data and examples}

We wish to write down the data required to construct such ``generalised irregular Betti spaces''.

First it is convenient to formalise the (well-known) notion of gluing 
quasi-Hamiltonian spaces. 
Given a quasi-Hamiltonian $G\times G\times H$-space $M$, we may fuse the two $G$ factors to obtain a quasi-Hamiltonian $G\times H$ space. Then we may reduce by the $G$ factor (at the identity of $G$) to obtain a quasi-Hamiltonian $H$-space, the {\em gluing} of the two $G$-factors.
Thus e.g. two quasi-Hamiltonian $G$-spaces $M_1, M_2$ may be glued to obtain a symplectic manifold (if it is a manifold) 
by gluing %
 their product:

$$M_1\glue{} M_2 = (M_1 \fus M_2)\spq G.$$

Now we will describe the required data.
Let $\cS$ be a finite set.
For each $s\in \cS$ choose a compact surface-with-boundary\footnote{In a more algebraic approach one should take the real oriented blow-up of a smooth projective curve, at some marked points.} $\Sigma_s$,
and a connected complex reductive group $G_s$.
Let $\cB_s$ be the set of boundary components of $\Sigma_s$
and let $\cB=\bigsqcup \cB_s$ be the set of all boundary circles.
(Write $G_b=G_s$ for any $b\in \cB_s$.)

For each $b\in\cB$ choose a semisimple element $A_b\in\g_b$ in the Lie algebra of $G_b$.
Let $H_b \subset G_b$ be the stabiliser of $A_b$.
Suppose we have chosen, for each boundary component $b\in B$, a product decomposition
$$H_b = \prod_{i\in I_b} H_i$$
of $H_b$ as a product of subgroups $H_i\subset H_b$. (We are not assuming that each $H_i$ does not decompose further, only that some decomposition has been chosen, possibly with only one factor $\# I_b=1$).
 
Let $I=\bigsqcup I_b$ be the disjoint union of the sets indexing the subgroups $H_i$. 
(Thus $I$ may be thought of as the set of boundary components obtained after gluing on all the pieces $\gah$ with  $G=G_b,H=H_b$. These components should be paired up, or glued to a conjugacy class, as follows.)

Choose a subset $K\subset I$ and a conjugacy class $\cC_i\subset H_i$ for each $i\in K$.

\begin{defn}
Given the data above, a gluing datum is an involution
$$\varphi:I\to I,\qquad \varphi^2 = 1$$
such that

1) $\varphi(i)= i$ if and only if $i\in K$, and

2) $H_{\varphi(i)} \cong H_i$ for all $i\in I$.
\end{defn}

Thus, given a gluing datum, the procedure to construct a symplectic moduli space is as follows.
Take the flat $G_s$-connections on each of the surfaces $\Sigma_s$ (with a framing at one point on each boundary component) and glue on a piece $\gah$ for each boundary component. This yields a quasi-Hamiltonian $H$-space,
where $H= \prod_{i\in I} H_i$.
Then glue together the factors $H_i$ and $H_{\varphi(i)}$ for all 
$i\in I\setminus K$, to obtain a quasi-Hamiltonian $\prod_{i\in K} H_i$-space.
This space is then reduced, at the value $\prod\cC_i$ of the moment map, to obtain a symplectic manifold (if it is a manifold). 

Note that we have not taken into account 
the choice of basepoints on the boundary circles, nor the choice of positive roots
appearing in the definition of $\gah$. These choices are important to understand the braid group action on the spaces (cf. \cite{bafi}) and may be encoded via Stokes representations of certain groupoids (similarly to \cite{smid} \S 3), but the point here is that  up to isomorphism the 
resulting complex symplectic manifolds will not depend on these choices. 
This is true for general reasons (``isomonodromy is a symplectic connection'' \cite{smid}) and follows in the present context from the fact that the two-form $\omega$ in Theorem \ref{thm: qh} does not depend on the semisimple element $A\in \g$.

\begin{eg}
Suppose $\# \cS=1$ so there is only one initial surface, and $K=I$, so there is no ``cross-gluing'', and $G_s=\GL_n(\IC)$ is a general linear group.
Then the above spaces are the Betti descriptions
(in the special case when the parabolic structures are trivial, and all irregular singularities have Poincar\'e rank one)  
of the hyperk\"ahler manifolds of \cite{wnabh}.
They are described in this way in \cite{saqh} in the case when each $A_b$ is regular (or zero).
\end{eg}

This example (and the usual case when each $A_b$ is zero) leads us to conjecture that such spaces are hyperk\"ahler in general.

\begin{rmk}

1) Note in general we have a semisimple Lie algebra element on {\em both} sides of each halo %
so the groups on each side of each halo only need to have isomorphic stabiliser groups.
Thus if we  ignore fission then the local picture at a halo is as follows:
$$M_1\glue{} \,_{G_1}\cA_H \glue{} \,_H\cA_{G_2} \glue{} M_2$$
where $M_i$ is a quasi-Hamiltonian $G_i$-space for $i=1,2$.

2) %
To recover the picture of the introduction suppose each set $I_b$ has just one element (so there is no fission), and also take one of the semisimple elements $A_b$ to be zero on one side of each halo.

3) Note that if $K$ is empty there will be no parabolic structures to worry about.

4) Note we have not ruled out the possibility of gluing two isomorphic factors $H_i\cong H_j$ of the same group $H_b$. 
\end{rmk}

\begin{figure}[ht]
	\centering
	\input{triangle.pstex_t}
\end{figure}

\begin{eg}
Here is an explicit example of the kind of generalised character varieties that arise.
A key feature is that more than one equation appears, in contrast both to the usual case and to the irregular case considered in \cite{saqh}.
Choose positive integers $a,b,c,g_1,g_2,g_3$. 
Set $n_1 = a+c, n_2=a+b, n_3=b+c$. Consider Riemann surfaces $\Sigma_i$ with genera $g_i$ for $i=1,2,3$ and each with one boundary circle. 
We will consider flat connections on $\Sigma_i$ with structure group $G_i :=\GL_{n_i}(\IC)$.
Consider ``block diagonal'' subgroups 
$$H_1 = \GL_a(\IC)\times\GL_c(\IC)\subset G_1$$
$$H_2 = \GL_a(\IC)\times\GL_b(\IC)\subset G_2$$
$$H_3 = \GL_b(\IC)\times\GL_c(\IC)\subset G_3$$
obtained as stabiliser subgroups of some semisimple Lie algebra elements with just two eigenvalues.
Denote the corresponding (triangular) unipotent groups as follows
$$
U_+, U_-\subset G_1,\qquad
V_+,  V_-\subset G_2,\qquad
W_+, W_-\subset G_3
$$
so $\dim(\U_\pm) = (n_1^2 - a^2 - c^2)/2 = ac$ and similarly $\dim(V_\pm)=ab, \dim(W_\pm)=bc$.
Denote by $\pi_a$ the projection $H_1\to \GL_a(\IC)$, and similarly for the other projections.
Thus $G_i$ is broken to $H_i$ for each $i$ and then we will glue the $\GL_a(\IC)$ factor of $H_1$ to that of $H_2$, and similarly for the others.
The resulting generalised character variety may then be described as follows (after some relabelling).
First on the three Riemann surfaces we have monodromy relations of the form:
\begin{align*}
[\al_1,\be_1]\cdots[\al_{g_1},\be_{g_1}] &= u_-h_1u_+\in G_1 \\
[\ga_1,\de_1]\cdots [\ga_{g_2},\de_{g_2}] &= v_-h_2v_+\in G_2 \\
[\eps_1,\ze_1]\cdots [\eps_{g_3},\ze_{g_3}] &= w_-h_3w_+\in G_3
\end{align*}
with $\al_i,\be_i\in G_1, u_\pm\in U_\pm, h_1\in H_1$ etc, where the square brackets denote the multiplicative commutator.
Then we have the gluing equations
\begin{align*}
\pi_a(h_1)\pi_a(h_2)&=\Id\in \GL_a(\IC),\\
\pi_b(h_2)\pi_b(h_3)&=\Id\in \GL_b(\IC),\\
\pi_c(h_1)\pi_c(h_3)&=\Id\in \GL_c(\IC).
\end{align*}
The resulting Betti space is the (affine GIT) quotient of this data:
$$
\left\{
(\al_i,\ldots,\ze_i, u_\pm,v_\pm, w_\pm, h_1,h_2,h_3)\st \text{monodromy and gluing equations}
\right\}\Bigl/ \Delta$$
where 
$$\Delta\cong \GL_a(\IC)\times \GL_b(\IC)\times \GL_c(\IC)\subset H_1\times H_2\times H_3$$
is the diagonal subgroup of $\Pi H_i = (\GL_a(\IC)\times \GL_b(\IC)\times \GL_c(\IC))^2$.
The group action is as follows: 
$H_1$ acts on $\al_i,\be_i,u_\pm, h_1$ by diagonal conjugation and fixes the other data (similarly for $H_2,H_3$) thus the subgroup $\Delta\subset \Pi H_i$ also acts and moreover preserves the six equations.
In other words the six equations define an affine subvariety of the affine variety
$$G_1^{2g_1}\times G_2^{2g_2}\times G_3^{2g_3}\times 
U_+\times U_-\times V_+\times V_-\times W_+\times W_-\times H_1\times H_2\times H_3$$
and we take the quotient of that by $\Delta$, i.e. we take the affine variety 
associated  to the ring of invariant functions.\footnote{as usual adding some sufficiently generic conjugacy classes ensures the resulting varieties are smooth (e.g. since in this way it is easy to ensure the data on each surface has no triangular decomposition, and so the action of the projectivisation of $\Pi H_i$ on the data satisfying just the monodromy relations, is free).}
Theorem \ref{thm: qh} and the general quasi-Hamiltonian yoga imply that (the smooth locus of) this variety is a complex symplectic manifold.

\end{eg}

\begin{eg}
Take two surfaces $\Sigma_1,\Sigma_2$ where $\Sigma_1$ is a two-holed sphere (with boundary components labelled by $0,1$) 
and with group $G=\GL_n(\IC)$ and with $A_0=0$ and $A_1$ regular semisimple, so that $H_0=G$ and $H_1=(\IC^*)^n$.
Take $\Sigma_2$ to be an $n$-holed sphere with group $\IC^*$ (and with arbitrary semisimple elements). 
Choose a (determinant one) conjugacy class $\cC\subset G$ to glue to the boundary $0$, and glue the $n$ boundary components of $\Sigma_2$ to the $n$-factors of $H_1$.
The resulting symplectic manifold is then (up to a covering) the symplectic leaf lying over $\cC$ of the Poisson Lie group $G^*$ dual to $G$ (cf. \cite{smapg} Lemma 4 and Proposition 23).
\end{eg}

\begin{rmk}[Irregular singularity viewpoint]\label{rmk: irreg}
Suppose we have a meromorphic connection on a $G$-bundle on a curve, where $G$ is a connected complex reductive group such as $\GL_n(\IC)$. To construct complex symplectic moduli spaces one fixes the formal isomorphism class of the connection at each singularity.
For simplicity suppose (similarly to \cite{wnabh}) that at each pole the connection is 
formally isomorphic to a connection of the form
\beq\label{eq: fnf}
A_\irr+\frac{\La}{z}dz
, \qquad\text{where}\qquad A_\irr = 
\left(\frac{A_k}{z^k}+\cdots+\frac{A_2}{z^2}\right)dz
\eeq
where the $A_i$ are commuting semisimple elements of the Lie algebra $\g$ of $G$ and $\La\in \g$ 
commutes with all the $A_i$
(thus we are assuming there is no ramification needed to obtain the normal form).
By definition the {\em formal monodromy} is $\exp(2\pi i \La)$;
this is the monodromy of the formal normal form \eqref{eq: fnf}. 
In the { regular singular} case we have $A_\irr=0$ and fixing $\Lambda$ 
amounts to fixing the local monodromy of the connection 
to be conjugate to $\exp(2\pi i \La)$ (the formal monodromy is conjugate to the actual local monodromy).
This can also be viewed in terms of a punctured Riemann surface or in turn more analytically in terms of a surface with boundary (either by removing a open disk or by performing the real oriented blow up at each singularity): we just consider 
nonsingular flat connections with fixed local monodromy conjugacy classes.
This yields an ``end'' of the surface.
Alternatively (not fixing $\Lambda$),
these surfaces with boundary may be glued together to obtain more complicated surfaces provided the monodromies match up.
More algebraically one may think of this heuristically as
gluing two curves (with regular singular connections) to form a node at the singular point, such that the residues $\Lambda$ add up to zero\footnote{
Going between these two viewpoints (nodal curves or sewing boundaries), is essentially the complex version of the exponentiation process of \cite{AMM}, together with the relation between quasi-Hamiltonian $G$-spaces and Hamiltonian loop group spaces.}.

Now in the irregular case we let $H\subset G$ be the subgroup stabilising $A_\irr$, so that $\La$ is constrained to be in the Lie algebra $\h$ of $H$ (and is only 
determined up to the adjoint action of $H$).
Thus the irregular part of the normal form yields a mechanism to reduce from $G$ to $H$. 
$\big(\!$
The formal monodromy is in $H$, and it is no longer necessarily conjugate to the actual local monodromy (in $G$), since a formal isomorphism to \eqref{eq: fnf} will not in general converge.$\big)$
First we may fix $\La$, thereby fixing the whole formal type, getting spaces as in \cite{wnabh, smid, saqh}, yielding what might be called an ``irregular end'' of the surface. 
Alternatively
if, as above rather than coming to an end, we now don't fix $\La$ 
we may try to sew surfaces with boundary (by taking the oriented real blow-up). 
The difference is that we should in the first instance glue on an $H$-connection with monodromy conjugate to $\exp(-2\pi i\La)$ ($\sim$ regular singular with residue $-\Lambda$), rather than a $G$-connection. Typically $H$ will decompose as a product and each factor may arise from a connection on a different surface.
More generally (and symmetrically) in the second instance we would glue on irregular connections with {\em formal} monodromy $\exp(-2\pi i\La)$.
The construction in the body of this article amounts to working out this idea in terms of Stokes data (i.e. the Betti realisation of such connections) in the case $k=2$, with the semisimple elements 
$A\in \g$ identified with the coefficient $A_2$ in $A_\irr$.
For larger $k$ this will be discussed elsewhere (the appendix of \cite{rsode} contains an additive/quiver analogue for general linear groups)---but let us mention that the gluing procedure 
looks to be precisely what is needed to build symplectic spaces of connections with multiple levels out of those with only one level.
\end{rmk}

\appendix
\section{Proof of Theorem \ref{thm: qh}}

First we will recall our notational conventions (largely from \cite{saqh}). %
$H,U_\pm \subset G$ are as in the body of the article and we denote the
corresponding Lie algebras $\h,\lu_\pm\subset\g$.
We have chosen a symmetric nondegenerate invariant
bilinear form $(\ ,\ ):\g\otimes\g\to \IC$. 
(Note that, since it is invariant,  $(\ ,\ )$ restricts to zero on
$\lu_\pm\otimes(\lu_\pm\oplus\h)$ 
and to a nondegenerate pairing on each of
$\lu_\pm\otimes\lu_\mp, \h\otimes\h$.)
The Maurer--Cartan forms on $G$ are denoted $\theta,\overline\theta\in\Omega^1(G,\g)$ 
respectively (so in any representation 
$\theta=g^{-1}dg, \overline\theta=(dg)g^{-1}$).
Generally if $\A,\cB,\cC\in\Omega^1(M,\g)$ are $\g$-valued
holomorphic one-forms on a complex manifold $M$ then
$(\A,\cB)\in\Omega^2(M)$ and 
$[\A,\cB]\in\Omega^2(M,\g)$ are defined
by wedging the form parts and pairing/bracketing the Lie algebra parts
(so e.g. $(A\al,B\be)= (A,B)\al\wedge\be$ for
$A,B\in \g, \al,\be\in\Omega^1(M)$).
Define $\A\A := \frac{1}{2}[\A,\A]\in \Omega^2(M,\g)$ (which works out
correctly in
any representation of $G$  using matrix multiplication).
Then one has 
$d\theta=-\theta^2, d\overline\theta = \overline\theta^2$.
Define $(\A\cB\cC) = (\A, [\cB,\cC])/2\in\Omega^3(M)$ (which is
invariant under all permutations of $\A, \cB, \cC$). The canonical 
bi-invariant three-form on $G$ is then
$\eta:= \frac{1}{6}(\theta^3)$.
The adjoint action of $G$ on $\g$ will be 
denoted $gXg^{-1}:=\Ad_gX$ for any
$X\in\g, g\in G$.
If $G$ acts on $M$, the fundamental vector field $v_X$ of $X\in\g$ is minus
the tangent to the flow 
$(v_X)_m = -\frac{d}{dt} (e^{Xt}\cdot m)\bigl\vert_{t=0}$, so that the map
$\g\to\Vect_M; X\to v_X$ is a Lie algebra homomorphism. 
(This sign convention differs from \cite{AMM} 
leading to sign changes in the quasi-Hamiltonian axioms and
the fusion and equivalence theorems.)

Recall that a complex manifold $M$ is a 
{\em complex quasi-Hamiltonian $G$-space}
if there is an action of $G$ on $M$, 
a $G$-equivariant map $\mu:M\to G$ (where $G$ acts on itself by
conjugation) and a $G$-invariant holomorphic two-form
$\omega\in \Omega^2(M)$ such that:

\noindent(QH1). 
$d\omega = \mu^*(\eta).$

\noindent(QH2).
For all $X\in \g$,
$\omega(v_X,\cdot\,) = \frac{1}{2}\mu^*(\theta+\overline\theta, X)
\in \Omega^1(M).$

\noindent(QH3).
For all $m\in M$: 
$$\ker(\omega_m)=\left\{(v_X)_m\ \bigl\vert\ 
X\in\g \text{ satisfies } gXg^{-1} = -X \text{ where } 
g:=\mu(m)\in G\right\}.$$

Recall that Theorem \ref{thm: qh}  claims that 
$\gah$ is a quasi-Hamiltonian $G\times H$ space.

\pfms (of Theorem \ref{thm: qh}).\ \ 
First we will check (QH1).
Write the moment map as $\mu = (\mu_G,\mu_H):M\to G\times H$, with 
$M=\gah$.
Since $\mu_G = C^{-1}pC$ we see $\mu_G^*(\theta)$ is conjugate to
$$\bar\ga + \cP -p^{-1}\bar\ga p$$
and so $\mu_G^*(\theta^3)$ equals $(\cP^3)$ plus three times
\beq\label{eq: coeff of 3}
(\bar\ga\cP^2) +(\bar\ga^2\cP)
- (\bar\ga p\bar\ga^2p^{-1})+(\bar\ga^2p\bar\ga p^{-1}) 
- (\bar\ga\bar\cP^2)+(\bar\ga^2\bar\cP) 
- 2(\bar\ga\cP p^{-1}\bar \ga p).
\eeq
On the other hand, from the definition of $\omega$ it is immediate that
$$2d\omega = \text{\eqref{eq: coeff of 3}} + 
d(\bar\cU_-,h\bar\cU_+h^{-1}).$$
Now if we expand both the last term here and $(\cP^3)$ (using the definition of $p$ to see $\cP = \cU_+ + u_+^{-1}\ch u_+ - p^{-1}\cU_-p$) we see
$$(\cP^3) = (\ch^3) + 3d(\bar\cU_-,h\bar\cU_+h^{-1}).$$
(Here $\ch = h^*(\theta_H)$ with $\theta_H$ the left-invariant Maurer--Cartan form on $H$.)
Thus since $\mu^{-1}_H(\theta_H^3) =- (\ch^3)$ we deduce (QH1), that
$\mu^*_G(\theta^3) + \mu_H^*(\theta_H^3) = 6d\omega$.

Next we will check (QH2). 
By linearity in $X$ we may check (QH2) separately for the actions of $G$ and $H$. First we consider just the $G$ action.
Choose $X\in\g$.
We will denote derivatives along $v_X$ by primes, 
so e.g. 
$\cP' = \langle v_X, \cP\rangle\in \Omega^0(M,\g)$
(and in any representation of $G$ we have
$\cP' = p^{-1}p'$ etc).
By definition $\mu_G^*((\theta + \bar\theta,X)) = 
(\mu_G^*\theta + \mu_G^*\bar\theta,X)$.
Expanding $\mu_G = C^{-1}pC$ we see
$$
\mu_G^*\theta + \mu_G^*\bar\theta
=
C^{-1}(\cP+\bar\cP)C + C^{-1}(p\bar\gamma p^{-1}-p^{-1}\bar\gamma p)C
\in \Omega^1(M,\g).$$

On the other hand considering the fundamental vector field of 
$X\in \g$ we see
$\gamma'=X, \bar\gamma' = CXC^{-1}, \cP'=\bar\cP'=\cU'_\pm=0.$
Thus
$$2\omega(v_X,\cdot) = 
(CXC^{-1},p\bar\gamma p^{-1}) - 
(pCXC^{-1}p^{-1},\bar\ga) + (CXC^{-1}, \cP+\bar\cP).$$
Clearly this rearranges into the above expression for
$(\mu_G^*\theta + \mu_G^*\bar\theta,X).$
Moving on to the $H$-action, first note that 
$\mu^*_H(\theta_H+\bar\theta_H) = 
-{\bar\ch}-\ch\in \Omega^1(M,\h).$
Then considering the action of $X\in \h$ we see
$\bar\gamma'=-X, \ga'=-C^{-1}XC, 
\bar\cU'_\pm = u_\pm Xu_\pm^{-1}-X, 
\cP'+\bar\cP' =  pXp^{-1} -p^{-1}Xp$
and in turn deduce:
 $$2\omega(v_X,\cdot)=
(u_-Xu_-^{-1} - X, h\cU_+h^{-1}) 
- (h^{-1}\cU_- h, u_+Xu_+^{-1} -X)
-(X,\cP+\bar\cP). $$
Upon expanding $\cP+\bar\cP$ (using $p=u_-^{-1}hu_+$)
this simplifies 
to  $-(X,\ch+\bar\ch)$ as required (by noting that $X\in \h$ pairs to zero with $\u_\pm$). This establishes (QH2).

For (QH3) it is now sufficient to check $\ker(\omega)\cap \ker(d\mu) = 0$ at each point (cf. \cite{ABM-purespinors} p.49). 
Thus choose a point $m\in M=\gah$ and a tangent vector
$X\in T_mM$. Suppose that $X\in\ker(\omega)\cap \ker(d\mu)$.
Since $X$ is in the kernel of $d\mu_H$ we have
$\ch'=0$ (here primes denote derivatives along $X$, so 
$\ch':=\langle h^*(\theta_H),X\rangle$).
Moreover $X$ being in the kernel of $d\mu_G$ amounts to the condition
$\bar\ga'+\cP'=p^{-1}\bar\ga' p$. Expanding $\cP'$ (and using $\ch'=0$) this becomes
\beq\label{star}
\bar\ga' +\cU'_+ = p^{-1}(\bar\ga'  + \cU'_-) p.
\eeq
Now we choose an arbitrary tangent vector $Y\in T_mM$ and denote 
derivatives along $Y$ by dots, 
so e.g. 
$\dot\cP = \langle Y, \cP_m \rangle\in \g$.
We then compute
\begin{align}\label{row1}
2\omega(X,Y) &= 
\left(
p^{-1}(\bar\ga'+\cU_-')p 
- p(\bar\ga'+\cU_+')p^{-1}
+\cU_-' - \cU_+' , \dot{\bar\ga}\right)\\\label{row2}
&+\left(
\bar\ga' + p^{-1}(\bar\ga'+\cU_-')p
, \dot\cU_+\right)\\\label{row3}
&-\left( \bar\ga'+  p(\bar\ga'+\cU_+')p^{-1}, \dot\cU_-\right)\\\label{row4}
&+\left( 
u_+\bar\ga'u_+^{-1} + h^{-1}u_-\bar\ga'u_-^{-1}h
, \dot\ch \right).
\end{align}
This should be zero for all $Y$; observe that each line is really an independent condition on $X$.
First \eqref{star} and its conjugate by $p$ imply that the right-hand side of \eqref{row1} is identically zero (we have already used $\ch'=0$ to compute \eqref{row1}).
Next \eqref{star} implies that \eqref{row2} equals
$\left(
2\bar\ga' + \cU_+' , \dot\cU_+
\right) = 2(\bar\ga', \dot\cU_+)$
so we see the $\u_-$ component of
$\bar\ga'$ is zero.
Similarly from \eqref{row3} we see the $\u_+$ component of
$\bar\ga'$ is zero.
Then \eqref{row4} implies the $\h$ component of $\bar\ga'$ is also zero, and so $\bar\ga' =0$.
Finally \eqref{star} now implies $\cU'_+=\cU'_-=0$, so, since all its components vanish, we see $X$ is indeed zero as required.
\epfms

\noindent{\bf Acknowledgments.}
Thanks are due to B. Kostant for kindly answering some basic questions concerning stabilisers of semisimple Lie algebra elements.
This research is partially supported by ANR grants 
08-BLAN-0317-01/02 (SEDIGA), 05-BLAN-0029-01 (GIMP).

\renewcommand{\baselinestretch}{1}              %
\normalsize
\bibliographystyle{amsplain}    \label{biby}
\bibliography{../thesis/syr} 

\def\cprime{$'$} \def\cprime{$'$} \def\cprime{$'$}
\providecommand{\bysame}{\leavevmode\hbox to3em{\hrulefill}\thinspace}
\providecommand{\MR}{\relax\ifhmode\unskip\space\fi MR }
\providecommand{\MRhref}[2]{%
  \href{http://www.ams.org/mathscinet-getitem?mr=#1}{#2}
}
\providecommand{\href}[2]{#2}
\begin{thebibliography}{10}

\bibitem{ABM-purespinors}
A.~Alekseev, H.~Bursztyn, and E.~Meinrenken, \emph{Pure spinors on {L}ie
  groups}, arXiv:0709.1452.

\bibitem{AMM}
A.~Alekseev, A.~Malkin, and E.~Meinrenken, \emph{Lie group valued moment maps},
  J. Differential Geom. \textbf{48} (1998), no.~3, 445--495, math.DG/9707021.

\bibitem{wnabh}
O.~Biquard and P.~P. Boalch, \emph{Wild non-abelian {H}odge theory on curves},
  Compositio Math. \textbf{140} (2004), no.~1, 179--204.

\bibitem{smapg}
P.~P. Boalch, \emph{{S}tokes matrices, {P}oisson {L}ie groups and {F}robenius
  manifolds}, Invent. math. \textbf{146} (2001), 479--506.

\bibitem{smid}
\bysame, \emph{{S}ymplectic manifolds and isomonodromic deformations}, Adv. in
  Math. \textbf{163} (2001), 137--205.

\bibitem{bafi}
\bysame, \emph{G-bundles, isomonodromy and quantum {W}eyl groups}, Int. Math.
  Res. Not. (2002), no.~22, 1129--1166, math.DG/0108152.

\bibitem{saqh}
\bysame, \emph{Quasi-{H}amiltonian geometry of meromorphic connections}, Duke
  Math. J. \textbf{139} (2007), no.~2, 369--405, math.DG/0203161.

\bibitem{rsode}
\bysame, \emph{Irregular connections and {K}ac--{M}oody root systems}, 2008,
  arXiv:math/0806.1050.

\bibitem{DMR-ci}
P.~Deligne, B.~Malgrange, and J.-P. Ramis, \emph{Singularit\'es
  irr\'eguli\`eres}, Documents Math\'ematiques, 5, Soci\'et\'e Math\'ematique
  de France, Paris, 2007.

\bibitem{Hit-sde}
N.~J. Hitchin, \emph{The self-duality equations on a {R}iemann surface}, Proc.
  London Math. Soc. \textbf{55} (1987), no.~3, 59--126.

\bibitem{MR91}
J.~Martinet and J.P. Ramis, \emph{Elementary acceleration and
  multisummability}, Ann. Inst. Henri Poincar\'e, Physique Th\'eorique
  \textbf{54} (1991), no.~4, 331--401.

\bibitem{Nak}
H.~Nakajima, \emph{Hyper-k\"ahler structures on moduli spaces of parabolic
  {H}iggs bundles on {R}iemann surfaces}, Moduli of vector bundles (Sanda 1994;
  Kyoto 1994), 1996, pp.~199--208.

\bibitem{Sim-hboncc}
C.~T. Simpson, \emph{Harmonic bundles on noncompact curves}, J. Am. Math. Soc.
  \textbf{3} (1990), 713--770.

\bibitem{simp-hfnc}
\bysame, \emph{The {H}odge filtration on nonabelian cohomology}, Algebraic
  geometry---{S}anta {C}ruz 1995, Proc. Sympos. Pure Math., vol.~62, AMS, 1997,
  pp. 217--281, math.AG/9604005.

\end{thebibliography}

\vspace{0.5cm}   
\'Ecole Normale Sup\'erieure et CNRS, 
45 rue d'Ulm, 
75005 Paris, 
 France

www.dma.ens.fr/$\sim$boalch

boalch@dma.ens.fr \qquad \qquad \qquad \qquad \qquad \quad\,\,  

\ 

\ 

\ 

\ 

\ 

[Journal Ref:  Ann. Inst. Fourier (Grenoble) 59 (2009) 7, 2669-2684]

\end{document}